
\magnification 1200  
\baselineskip=12pt
\input amstex
\documentstyle{amsppt}
\document
\def\ju{\vskip .8truecm plus .1truecm minus .1truecm}

\def\sju{\vskip .4truecm plus .1truecm minus .1truecm}
\def\qq{\lq\lq}
\def\lqq{\qq}
\def\pr{^\prime}

\def\lra{\longrightarrow}
\def\se{\subseteq}
\def\lan{\langle}
\def\ra{\rangle}
\def\RR{\Bbb R}
\def\prl{\partial}


\NoBlackBoxes
\vsize 22truecm
\font\sc =cmcsc10
\def\noi{\noindent}

\vskip 0.3truecm
\hfill February 12 2013
\sju
\centerline{ \sc \bf Feedback from students' tests as a tool in teaching }
\centerline{\bf A case study%
    \footnote  {The first version of this paper was presented at the joint American Mathematical Society/Mathematical Association of America Meeting, Boston, January 4, 2012. The paper was presented at the 12th International Congress on Mathematics Education (ICME12), Seoul, South Korea, July 8-15, 2012.} }
\vskip 0.5truecm
\centerline{\it Radoslav M.  Dimitri\'c}
\centerline{dimitricr\@member.ams.org}

\vskip 1.3truecm

{\bf Introduction.} A close look at students' 
written work on examinations offers a wealth of information about their performance, their knowledge of the subject, their strengths, weaknesses and misconceptions, and their overall level of mathematical skills and abilities. This information can be  used to better ascertain how mathematical  concepts taught in class were understood by students and to suggest  various approaches that could be used to improve teaching. The process of analyzing students' work, although time consuming, can be rewarding and can have positive results for both students and instructors.

My aim is fourfold: (a) To give a rationale for analysis of students' work, especially tests (b) to give possible explanations of students' performance; (c)  to revisit aspects of teaching some mathematical constructs, such as function and, (d) to provide possible ways of improving the teaching of concepts considered in this paper.

\noi $----------$
\sju

{\it Mathematics/Education Subject Classification} (MSC 2010): Primary: 97C70, 97D40, 97D70, 97I40; Secondary: 97D66, 97D80, 97I40 / (MESC 1999) Primary: B40, B50, C70, D40,  D70,  Secondary: D60, D80, E40, I20, I40, I60


\vfill
\eject

In this paper, I present the core of my analysis of one final examination in a Calculus III course for science majors (by way of answers about their majors/career goals, 42.8\% of the students were chemistry or biochemistry majors, 14.3\% math and math ed majors, 10.7\% of computer science majors, 10.7\%  mechanical or civil engineering, 7.1\% medicine or biomedicine; remaining students were  
pharmacy, history teaching, biostatistics, undecided and they amounted  amounted to 14.4\%.  
Judging by a diagnostic test given at the first class meeting and made specifically for this course (the diagnostic is provided in the appendix) these science students were not mathematically advanced, and showed considerable weakness on many points. Diagnostic test score of about 90\% would indicate that the student is prepared for the course and the score below 60\% requires instructor's serious discussion as to what to do about the required material. In this light the mean score of 35\% (median 38\% and standard deviation %
As is often the case, there were difficulties related to the concept of functions, their properties and their operations. 
There is a considerable corpus of literature that addresses the issue of teaching and more often students' understanding of the notion of function (a sample is Breidenbach et al (1992), Dubinsky \& Harel (1992), Horscovics (1982),  Oehrtman et al (2008), Piaget et al (1977), Sfard (1987), Tall, D. \& Vinner, S. (1981), Thomson (1994), Vinner (1989, 1992)) and a survey in Leinhardt et al (1990).
The fundamental starting points in their research is in the field of psychology. 

One strand of that research concentrates on students' perceptions (or subconscious metaphors) of the notion of function, dividing it into the process conception vs. action conception categories. After Dubinsky \& Harel (1992, p.85), \lqq an {\it action conception} of function involves the ability to plug numbers into an algebraic expression and calculate. It is a static conception in that the subject will tend to think about it one step at a time (e.g. one evaluation of an expression). A student whose function conception is limited to actions might be able to form the composition of two functions, defined by algebraic expressions, by replacing each occurrence of the variable in one expression by the other expression and then simplifying; however, the students would probably be unable to compose two functions that are defined by tables or graphs." In contrast, \lqq a {\it process conception} of function involves a dynamic transformation of quantities according to some repeatable means that, given the same original quantity, will always produce the same transformed quantity. The subject is able to think about the transformation as a complete activity beginning with objects of some kind, doing something to these objects, and obtaining new objects as a result of what was done. When the subject has a process conception, he or she will be able, for example, to combine it with other processes, or even reverse it. Notions such as 1 -- 1 or onto become more accessible as the student's process conception strengthens." 

My own work with students certainly supports importance of perceptions and metaphors both in learning and teaching, however perhaps more fundamental issues need to be addressed prior to or simultaneously with the problem of perceptions and those involve focus on teaching and learning precise \lqq formal" definitions of these concepts, that in turn, depend on precise, formal definitions of notions of two and higher dimensions (and representations in coordinate systems). Perhaps a notion close to this aspect is \lqq concept definition" as discussed in Tall \& Vinner (1981), Vinner (1991), Vinner \& Dreyfus (1989). However these works concentrate primarily on \lqq concept image" which is undoubtedly important, however the metaphors (such as in \lqq concept images") cannot be made operational and stable, unless there is definite command of relevant and definite definitions these metaphors should map. By metaphor, I understand the classical (and somewhat modernized) definition of metaphor originating with Aristotle in his {\it On the Art of Poetry} and {\it Rhetoric} \lqq carrying of expression from one object onto the other..." In short, metaphor consists of three parts: the first object, the \lqq carrying" i.e the process of connecting the two objects,  and the second object." Metaphors are different if at least one of their three parts differs from the corresponding part of another metaphor.

The present paper gives an example of analysis of students tests in a default classroom setting. Among other things, I make suggestions as a practical guide in teaching the subjects involved, but I also point out areas worthy of further education research. Notable studies exist with related themes and I mention here Selden et al (1989, 1994, 2000) where the authors test (outside of the classroom setting) self-selected sample of students (after they complete the usual calculus sequence) on routine and non-routine problems. 

\ju
\centerline{\bf The test, students' solutions and analysis}
\ju
Seven  concrete test problems are presented here (to be completed in 100 minutes -- the time that is perhaps tight to complete all the problems). By departmental decision, graphing calculators are allowed, starting from the most elementary remedial courses to post calculus courses. My comments accompany each individual problem with final conclusions at the end of the paper. Students' solutions are quoted verbatim, without grammatical corrections. 

When making the exams, I try to adhere to some extent (and to the extent possible in classes at this level) to Bloom's taxonomy (Bloom (1956), and its revision, Anderson et al (2001)) and devise exam questions of different difficulty to test various aspects of students' knowledge. The issue of the difficulty level and all the ramifications is a complex issue. I expound on this topic in more detail in Dimitric 2012b, while adopting here what is possibly the crudest (but still useful) graduation in Bloom's two-dimensional taxonomy array. The first two problems were meant to be of the simplest kind and students' scores on the problems show that my gage of relative difficulty was correct. This however does not mean that students' performance was appropriately good, as much as I expected or hoped for -- the mean scores on the individual problems, out of maximum 15,  were as follows (Problems 1 to 7, in order): 8.04, 6.39, 2.96, 3.52, 0.96, 2.22, 3.39 while the respective medians were: 8, 7, 1, 3, 0, 1, 3. My rating of problem 5 was not the highest difficulty, for I deemed it not so time consuming, however, students' performance on this problem makes it the hardest problem. The reasons for this are manifold, the main being the fact that the continuity concept (and function in general)  is an orphan in calculus education and another reason is that solving this problem required a bit of reasoning in steps. While these problems were not out of the ordinary in comparison to problems worked out in class or assigned for homework, the problems may be more difficult than what the students have seen in previous courses and or courses taken with other instructors. We learn from Selden et al (2000) that what they term \lqq non-routine" calculus I problems pose considerable difficulty even for students who completed the calculus sequence (including differential equations). Perhaps peculiarly, these authors report that, in place of attempting to solve the problems by calculus techniques,  the students use algebra,\footnote   {The issue of use of more pre-calculus in calculus was addressed in Dimitric (2001); that approach however requires a greater sophistication than is exhibited from a typical calculus student} or guessing, etc in their attempts. As it is, the problems are halting already at Bloom stage one complexity level\footnote   {Selden et al report that their routine stage is well fulfilled by the students they test, however their testing was on the material several courses earlier than the course levels the tested students were at, at the time of testing} and this is the area that needs much more attention in research literature. 

Since this kind of analysis is also extended to midterm exams (for any given class), the information I obtain is used to improve teaching in a number of ways. 

\sju

{\bf Problem 1.}  Given three points in 3D space: $A(2,0,a), B(0,1,a)$, $C(-1-1,0)$, where $a$ is an arbitrary given constant. 

a) Calculate the vector/cross product $\vec{AB}\times \vec{AC}$.

b) Calculate the scalar/dot product $\vec{AB}\cdot \vec{AC}$

c) For what values of $a$ are vectors $\vec{AB}$ and $\vec {AC}$
   (i) perpendicular, (ii) parallel?

d) Write an equation of the plane that contains points A,B,C. 
\vskip 1truecm

Normally, analytic geometry is not a part of the third calculus; rather, it is often \lqq reviewed" in previous calculus courses, and at any rate is included in most \lqq coffee table" calculus textbooks. My department decided to include it in calculus III as a \lqq review" and rightly so, given the level of mathematical skill of students taking such courses. There is however a drawback in this kind of remediation since the topics that are left for the end of the course are sometimes, if not half of the time (depending on the class), left out for sheer lack of time. But those topics are (or should be) the core of calculus III -- Green's, Stokes' and divergence theorems. This course is fairly standard in the US -- it includes topic on vector-valued functions and surfaces, partial derivatives, optimization and Lagrange multipliers, double and triple integrals and their properties and techniques of evaluation, vector fields, line integrals, surface integrals that all culminate in Green's, Stokes' and Divergence theorems. 

This problem was meant to test the elementary 3D concepts at the level of definitions and simple steps. I have used a constant $a$ in order to prevent \lqq button pushing" on programmable calculators. This causes some problems to some students who still think that a constant must be given as a concrete decimal number (most often an integer), not a letter, which (according to many students) is an unknown variable, by default. This kind of problem could be eliminated if there were beginning courses that would teach basic mathematical culture, including the concepts of a constant and variable; in fact such courses should be made mandatory, at least for mathematics majors. 

A number of students stumble upon an objective problem, namely the fact that we represent vectors both by their end-points (that are each ordered triples) as well as by a single triple of vector's coordinates (that are the difference of those end-points). This is an objective confusion since localization and affinity  of 3D space make this \lqq gluing" possible. The instructor should thus spend some time in elaborating the concept of vector translation which produces the same vector (strictly an equivalent vector) that begins at the origin. And judging by further errors made, the instructor should also make sure that students remember that you need to subtract the initial point from the end point of the  vector to get the vector's 3D components: 
$\vec{AB}=(-2,1,0)  $ and $\vec{AC}=(-3,-1,-a)$. 

One idea to avoid this confusion is perhaps to use round brackets for ordered triples representing points and to use angle brackets to represent vectors, however I did not find that this alone can remedy the problem.  

The errors even at this first step are in simple arithmetic (subtracting greater number from smaller), mishandling of minus signs and parentheses, multiplications by 0 treated as multiplications by 1, and getting vectors to be the opposite of what they should be. Arithmetic errors come somewhat from the exam pressure and to a greater extent because of addiction to calculators that result in students being unable to do even the simplest arithmetic operations in their heads. I have mentioned elsewhere [Dimitric 2003, 2009] that students did a bit better when they were not allowed calculators into the exams.

 One student \lqq finds" the vectors by dot multiplying their endpoints, which is in the domain of confounding (with the dot product).
I wanted to check that they know definitions of the cross (vector) and the dot (scalar) products and how to compute them if vectors are given in the component form. In the scalar product all that is required is component-wise multiplication and addition, but the vector product involves a bit of a twist in that students have to either use the appropriate determinant or express the vectors through the basis vectors $i, j, k$ and use properties of cross multiplication. A majority of the students used the determinant and perhaps it may be due to the fact that I used the determinant more often than other methods in computing the product when working out exercises at the blackboard. That should then imply that including more computations with vectors spanned by the basis may be helpful. 

In calculating the cross product determinant (which computes to be $-ai-2aj+5k$), the usual arithmetic errors occur, now including calculation of the three by three determinant. Not all students remember the determinant definition of the cross product. The student that \lqq found" vectors by multiplying their endpoints also formed a 2 by 3 \lqq determinant" and did some inventive multiplications inside to get his incorrect result. 

Not all students remember the coordinate definition (computation) of the scalar product (equal to 5 in this case). Here is how one student does it: 
$(-2,1,0)\cdot (-3,-1,-a)=(6,-1,0)$. 

Part c) was meant to check whether students learned the simple connection between the two products with vectors and orthogonality and parallelism. Some confuse the two and relate parallel to  $\vec{AB}\cdot \vec {AC}=0$ and perpendicular to $\vec{AB}\times \vec {AC}=0$ (many may miss the fact that the latter $0$ is a vector). Here the instructor has more work to do in devoting more time in explaining the distinctions between the two products, such as that one is a scalar and the other is a vector and spend more time on the geometric definitions of the two (in addition to their consequent computational definitions). One may add mnemonic rules to aid the students to remember which is which and how they are related to these important geometric conditions. 
Another student confuses perpendicularity with the condition for perpendicularity of two lines; he states: \lqq parallel when $a=0$.  Can't be perpendicular because $a$ cant [sic] become a negative reciprocal..."  That $a$ should be $0$ for parallel or perpendicular is frequent. One student computes the dot product correctly to be 5 and states correctly the condition of perpendicularity, but claims that the two vectors are always perpendicular; this could be simply a lack of concentration on the student's part? Another student \lqq solves" this as follows: 
\lqq $\vec{AB}$ and $\vec{AC}$ parallel if they point in same or opposite direction, since $\vec{AB}=\lan -2,-1,0\ra $ and
$\vec{AC}=\lan -3,-1,-a\ra $, $\vec{AB}$ and $\vec{AC}$ are not parallel because first component is unequal."  Perhaps one needs to spend more time on the concept of parallelism by way of one vector being a scalar multiple of the other, resulting in all corresponding coordinates being proportional. 

There were also students who overcomplicated, at least with some of their solutions. Calculus II or III students often take linear algebra concurrently and are only too eager to demonstrate what they are exposed to in that class, Gauss elimination being one of the favorites, since considerable attention is paid to it in linear algebra courses and it is a fairly algorithmic procedure that students always favor. Here is how one student approaches questions in part c) (decomposition into perpendicular and parallel components of a vector is also in his mind, for it may have been a part of a review session before the exam): 
 i) $v_{\perp}+v_{||}=||v||$, $a=-5$, ii) $ a=2.25$. To get this $a$, he does Gauss elimination in the matrix made up of rows $\vec{AB}$ and $\vec{AC}$ and gets matrix with rows $(1,-1/2,0)$ (incorrect) and $(0,-5/2,-a)$ and then matrix with rows $(1,-1/2,0)$ and $(0,1,2a/5)$ and then rows $(1,0,9a/10)$ and $(0,1,2a/5)$, then equates $9a/10=2a/5$ and  gets $45a=20a$ from which he derives $a=2.25$. Ending up with an incorrect answer may be prevented, with some chance of success, by the instructor telling students to check their answers by plugging in the values they get. 

Part d) was meant to give a 2D/3D analogue of the point-slope equation in the Cartesian plane. The plane is $ax+2ax-5z+3a=0$. Although students may have incorrectly found their cross product, if they wrote the correct form for the equation of the plane, they received most of the credit for this part.  Some confused the plane equation with that for a line through a point. Thus, this first step of simply memorizing the plane equation is worth working on. Certainly, memorization comes with a working out a good quantity of exercises that utilize the equation or concept to memorize; in addition my advice to students was to make an index card for a notion they are learning and write definitions and basic properties on one side and a typical example(s) using the notions on the other side of the card.  
\ju
\ju


{\bf Problem 2.} Find a function $f(x,y)$ that satisfies the equation $\nabla f=(y^2, x) $ or prove that such a function does not exist. (You may assume that the function has partial derivatives of any order.)
\sju

This exercise checks basic properties related to partial derivatives  and also whether students know the conditions, namely that $\prl F_1/\prl y=\prl F_2/\prl x$ and that the components should have continuous partial derivatives. If these conditions are satisfied, then they should also demonstrate whether they can find such a potential function $f$. To simplify considerations on smoothness of partial derivatives, I made the assumption on the existence of all the partial derivatives which then imply their continuity. If such an $f$ exists, then, we have to have 

\centerline{$\nabla f=(\prl f/\prl x, \prl f/\prl y)=(y^2,x)$, i.e. $\prl f/\prl x=y^2$ and $\prl f/\prl y=x$,} 

\noi consequently $\prl^2f/\prl y\prl x=2y$, $\prl^2 f/\prl x\prl y=1$. Since we are assured that all partial derivatives exist and are therefore continuous, a Clairaut's theorem asserts that these two mixed second order partials must be equal (which does not happen here, apart from  the points on a single horizontal line $2y=1$).

Students also had a ready-made result that was used in a few exercises in class, namely that if $(F_1,F_2)=\nabla f$ is a gradient vector field (and $f$ is a potential function) then one must have $\prl F_1/\prl y=\prl F_2/\prl x$, which would again lead to the same conclusion of no solution. 

However, the mathematical results on conditions and existence are not always welcome or easily adopted by science students; they prefer to \lqq do something" -- compute a number, a value of a function or such (which may be classified as an \lqq action conception of function", Dubinsky \& Harel (1992)). Thus, some students can ignore conditions and go straight to solving the system: (1) $ \prl f/\prl x=y^2, \prl f/\prl y=x$  (2). After integrating the first equation  with respect to $x$ one gets $f(x,y)=y^2x+\phi(y)$, hence $\prl f/\prl y=2yx+\phi\pr(y)$, which again leads to non-existence of the solution ($2y=1$). 

Indeed students' solutions (or incorrect solutions) show that they mostly attempted to find the solution by solving the system of the two pde, by finding (in)appropriate anti-derivatives. Thus students get that the function should be $xy^2+C$ after integrating the first equation with respect to $x$, likewise, after integrating the second equation with respect to $y$,  that the potential function  should be equal to $xy+C$, without realizing that these $C$'s are functions of, respectively, $y$ and $x$, rather than an identical constant for both (for some the $C$ is not there altogether). They still get the correct answer claiming that these two formulas are not the same. A couple of the students had a different \lqq idea" namely that $f$ is the sum (or the difference) of the two partial anti-derivatives that they obtain; $f(x,y)=xy^2\pm xy$. This error is due  most likely to shortness of time students have to separate and digest different topics they are learning at a fast pace for them. They have been exposed to various results with sums of partial derivatives, such as in the multy-variable chain rule, and they likely have that confusion in this regard. On the other hand firm grasp of definition of a function, and its domain and codomain and integration and differentiation with respect to one variable would have prevented more blatant errors. Thus keeping in mind that $f:\RR^2\lra \RR$, i.e. that the domain of $f$ is the (two-dimensional) plane would already be a good medicine.

A work due for the instructor is to emphasize repeatedly that when integrating with respect to one variable, the other variables are considered to be constants, but the \lqq constants" that come out of such integrations are only constant with respect to the variable of integration, and, that they are in fact functions of the remaining set of variables. That is not unlike introduction of partial derivatives. Generally the transition from one dimensional to two dimensional case deserves much attention and the instructor should not spare time in making sure that students get that very basic, but important understanding. Even with one variable, the anti-derivative constant $C$ is usually relegated to either not writing it down or to routine addition to the anti-derivative that they find, yet one should spend a bit of time on the constant as well to relate to students that an anti-derivative is in fact a family of functions that differ by constants and that $C$ may take different values, although we use that one letter for its label (in case of discontinuities the constant may be defined piecewise).

Another effort the instructor should make is to advertise conditions for these results as an essential part not only of the theory they learn, but as an essential part of their problem-solving skills.

A rare student who appealed to Clairaut's theorem forgets the crucial condition, namely the assumption on the continuity of the second partials, in order to make a conclusion about their equality. That, despite my efforts  to preempt this mistake while teaching them the topic, by pointing out the importance of the continuity condition; textbooks usually place an example (often in exercises section) that illustrates that continuity condition is needed. 

One can also see that there are a few (luckily only a few) students who seem not to  have learned yet what $\nabla$ means. That can be seen when one student  integrates the first equation with respect to $y$ and the second with respect to $x$. Another error one student makes (after finding the anti-derivatives) is that he writes $f(x,y)=(xy^2,xy)$ thus confusing the roles of the players in the question. Knowing that the codomain of $f$ is one-dimensional would have prevented that error and that, yet again would have been clear if the student wrote all three parts of the function: $f:\RR^2\lra\RR$.
\ju
\ju


{\bf Problem 3} Express partial derivatives $\prl f/\prl r, \prl f/\prl \theta$ of a function $f(x,y,z)$ in terms of $\prl f/\prl x, \prl f/\prl y, \prl f/\prl z$ where $(r,\theta, z)$ are cylindrical coordinates. 
\sju

The purpose of this exercise is to test students on the routine of the chain rule in several variables (three), each of which depends on three variables. Instruction was geared towards showing that multi-variable chain rule is actually a natural extension of the one variable case and that, for that matter, the chain rule formulas 

$$\frac{\prl f}{\prl r}=\frac{\prl f}{\prl x}\frac{\prl x}{\prl r} + \frac{\prl f}{\prl y}\frac{\prl y}{\prl r}   + \frac{\prl f}{\prl z}\frac{\prl z}{\prl r}      \eqno (1)$$

$$\frac{\prl f}{\prl \theta}=\frac{\prl f}{\prl x}\frac{\prl x}{\prl \theta} + \frac{\prl f}{\prl y}\frac{\prl y}{\prl \theta}   + \frac{\prl f}{\prl z}\frac{\prl z}{\prl \theta}      \eqno (2)$$
are understandable in view of the one variable chain rule. Just writing down these formulas earned them 1/3 of the points for this problem. They also had to know the cylindrical coordinate formulas, from which then they had to find the needed partial derivatives 
$$ \frac{\prl x}{\prl r}=\cos\theta, \frac{\prl y}{\prl r}=\sin\theta, \frac{\prl z}{\prl r}=0, \frac{\prl x}{\prl \theta}            =-r\sin\theta,  \frac{\prl y}{\prl \theta}=r\cos\theta, \frac{\prl z}{\prl \theta}=0 
$$
(this would earn the student almost another third of the points). They would be ready to substitute these partials into (1) and (2) to get the final answer:

$$
\frac{\prl f}{\prl r}=\frac{\prl f}{\prl x}\cos\theta+\frac{\prl f}{\prl y}\sin\theta; 
\frac{\prl f}{\prl \theta}=\frac{\prl f}{\prl x}(-r\sin\theta)+\frac{\prl f}{\prl y}r\cos\theta.
$$
\vskip 0.7truecm

A number of errors appear here. A few students in fact do not quite know the cylindrical coordinate formulas yet (Bloom's step one). Confusions are with the polar coordinates in two dimensions $x=r\cos\theta, y=r\sin\theta$ (omitting the third equation $z=z$, or writing $z=r$); this is also confused with the equation of the right circular cylinder $x^2+y^2=r^2$ or even $=z^2$, since the thinking is that $z=r$. At least one student confused cylindrical coordinates with spherical ones, undoubtedly attributed to the fact that, in a relatively short time, both of these coordinate sets were used in various problems and topics covered in class. 

Getting  incorrect derivative from cylindrical coordinates: $\partial x/\partial r=\sin\theta$ and $\partial y/\partial r=-\cos\theta$, $\partial z/\partial r=1$ was also observed. 

Some have not learned the chain rule formulas yet, even though they know the cylindrical coordinates formulas. Thus, students \lqq invent" their own versions of what they have not yet assimilated; for instance: $\partial f/\partial x=-r\sin\theta, \partial f/\partial y=r\cos\theta$, $\partial f/\partial z=0$. Or  $\prl f/\prl r\cdot  \prl r/\prl x=\prl f/\prl x$ and $\prl f/\prl \theta\cdot \prl \theta/\prl y=\prl f/\prl y$ and $\prl f/\prl z=\prl f/\prl z$ no change.  And $\prl x/\prl\theta \cdot \prl\theta/\prl r=\prl x/\prl r$ and $\prl y/\prl\theta\cdot  \prl\theta/\prl r=\prl y/\prl r$. 

Here instructor should spend much time in elaborating on function/variable dependencies of the various functions and variables involved, in fact making sure that students feel comfortable with functions and their compositions. The latter is especially the case in point since a good number of students are still struggling with ordinary, one-variable compositions of functions (and one variable chain rule) and compositions of functions of several variables require special attention, even if they are a very natural extension of the one-variable progenitor. One simple condition, namely that the composition $g\circ f$ of two functions exists only if $Range\,\, f\se Domain\,\, g$ is entirely missing from calculus textbooks, mainly because of sole reliance on \lqq intuitive" notions about functions (be it action, process, or \lqq covariation"). The objective problem is in that finding the range of a function is difficult, but this is exactly where more work should be invested, in addition to introduction of a notion more relaxed than the range, namely that of codomain.  The lack of understanding of these compositions is visible from the following incorrect reasoning: 
\lqq Since we are expressing the partial derivatives respectively in terms of $f(x,y,z)$ to $(r,\theta,z)$ that must mean that  $\prl f/\prl r=\prl f/\prl x$, $\prl f/\prl\theta=\prl f/\prl y$, $\prl f/\prl z=\prl f/\prl z$ respectively and the function $f(x,y,z)=f(r,\theta,z)$." 

Certainly it is a \lqq stopper" if they do not know formulas (1) and (2) (Bloom's level 1). Again, making a point that these formulas are an extension of the one-dimensional case would likely ease their remembering the multi-variable chain rule formulas. 

An important subtlety is our notation, namely one writes $f(x,y,z)$ and, after substituting $x,y,z$ as functions of other variables $r,\theta,z$, we may erroneously use the same symbol and write $f(r,\theta,z)$, an infraction that could be remedied when a seasoned mathematician is doing mathematics by knowing that it is an abuse of notation. However students will not pick up on that subtlety and they will instantly conclude that $f(x,y,z)=f(r,\theta,z)$ and thus $(x,y,z)=(r,\theta,z)$. At least one student must have concluded this since he states that  $\prl f/\prl r=\prl f/\prl x$ and $\prl f/\prl \theta=\prl f/\prl y$ and $\prl f/\prl z=\prl f/\prl z$, as mentioned above.

There were some students who displayed the chain rule equations, and the cylindrical coordinates but were nonetheless unable to complete the appropriate partials from the cylindrical coordinates. But then, there was at least one student who knew the cylindrical coordinates, knew how to compute the appropriate partials from them, but did not know what to do with it, since he was missing basic chain rule equations (1) and (2). One student is confounding the chain rule in this problem with the substitution formula (change of variables), thus he writes $\iiint f(r,\theta,z) dzdrd\theta$. 

 Here doing chain rule drills with a good number of examples seem to be the only way to learn the multy-variable chain rule. 

One characteristic problem with this exercise must have been, for most students, the fact that there was no concrete formula for the function $f$, which made students attempt to invent it, arguably, because of the action conception of function, but even with the formula, knowledge of the (dimensions of) domain and codomain would have had helped greatly). The students were led astray also by the phonetics they had here. Thus \lqq cylinder" and \lqq cylindrical" induces one student to  write a circular paraboloid equation $x^2+y^2=z$ (she likely wanted the equation of a  cylinder) and then she  writes $f(x,y,z)=x^2+y^2-z$, $f(r,\theta,z)=(r\cos\theta)^2+(r\sin\theta)^2-z=2r^2-z$ and also  $\prl f/\prl x=2x$, $\prl x/\prl r=\cos\theta$, $\prl f/\prl y=2y$, $\prl y/\prl\theta=r(\cos\theta)$ then
$\prl f/\prl r=\prl f/\prl x\cdot\prl x/\prl r=2x\cos\theta=2x\cos\theta$,  $\prl f/\prl\theta=\prl f/\prl y \prl y/\prl\theta =2y r\cos\theta=2yr\cos\theta$, thus inventing the \lqq chain rule" (notwithstanding algebraic errors).  Another student writes equations of cylinders or spheres (radii $=r$) as the formula for $f(x,y,z)$ (which is not given), again confounding cylindrical coordinates with inventing a formula for $f$. Yet another student writes $f(x,y,z)=x^2+y^2+z$ $f_x=2x, f_y=2y, f_z=1$ $r=\sqrt{x^2+y^2}, \theta=\tan^{-1}(y/x), z=z$.
\ju
\ju


{\bf Problem 4} Given points A(0,1), B(1,2) and C(3,0) (in 2D) and the curve $\gamma$ that consists of the portion of the parabola  $y=1+x^2$ between points A and B and the straight line segment connecting B and C,  find the line integral 
$$\int_\gamma\vec\bold F\,\cdot d\vec\bold s$$
where $\vec F(x,y)=(x^2-y, y^2+x)$ and the curve is traversed in the direction from A to B to C. 
What would be the value of the integral if the orientation of the curve $\gamma$ were from C to B to A?

\sju

This one is meant to test whether students have gotten basics of the line integral of a vector field. Parts of the solution(s) that were assigned partial credit were: Drawing the graph of the path and recognizing that the integral over $\gamma$ is the sum of the two pieces -- over the parabolic path and over the line path: $\int_\gamma=\int_{\gamma_1}+\int_{\gamma_2}$ (2 points). Then the students should know how to parametrize the two paths, for instance: $\gamma_1: x(t)=t, y(t)=1+t^2, t\in[0,1]$ and $\gamma_2: x(t)=t, y(t)=-t+3, t\in[1,3]$. One of course needs to know the basic technical integration formula $\int_\gamma \vec\bold F(x,y)\cdot d\vec\bold s=\int_a^b (\vec\bold F(x(t),y(t))\cdot (x\pr(t),y\pr(t))dt$ and to find correct derivatives $x\pr=1, y\pr=2t$, for $\gamma_1$ and $x\pr=1, y\pr=-1$ for $\gamma_2$ (3 points). Routine polynomial definite integration will then produce $\int_{\gamma_1}\vec\bold Fd\vec\bold s=2$ and $\int_{\gamma_2}\vec\bold Fd\vec\bold s=0$ (6 points), so that $\int_\gamma=2+0=2$ (2 points). The orientation question tests whether a student knows that reversing the path orientation changes the sign of the value of the integral i.e. $\int_{\gamma^-}=-\int_{\gamma^+}$ (2 points). 

Apart from the usual arithmetic and algebraic errors, some students are unable to sketch the (entire) integration path (one student attempts to draw it in 3D, undoubtedly influenced by a good number of 3D images of surfaces and curves done on other topics). Ability to draw these simple curves should be a prerequisite for taking any of the calculus courses (see here Dimitric (2012a)). Again, some additional time should be spent on functions defining curves and surfaces in different dimensions; it is highly recommended that functions are presented not only as formulas, but as ordered triples consisting of domain, codomain, and the formula(a rule of correspondence). Analysis of the dimensions of the domain and codomain would then aid  the student in recognizing the context as to whether he is dealing with a surface or a curve and whether it is two, three, or more dimensions. Striving to make students switch from action to process conception of function would not help, unless students actually were well familiar with the complete functions, namely domain, assignment rules and codomain, or luckily the range. 

There were a couple of students who only drew the graph of the path and did nothing else. 
When it comes to parametrizations the mistakes are in not parametrizing the line or the parabola piece, or not determining correct limits of the parameter in either of the two pieces. Some students however see no need for parametrizations and use the Cartesian coordinates of the two paths (and some may  have the line segment Cartesian equation incorrectly). This is because those students likely confound this with the double integral exercises they have been exposed to. Problems with parametrizations should be diminished by the instructor's explanation why parametrizations are needed, as well as working out some simple but also some intricate parametrizations, while elaborating on the conditions, advantages or disadvantages of one parametrization over the other. 

A number of students then use an incorrect formula to integrate this vector field, namely $\int_\gamma \vec F\cdot d\vec s=\int F(c(t))\cdot ||c\pr(t)|| dt$, which is a formula for the line integral of a scalar function $F$. Yet again, dimensional analysis here is of great help to decide the correct context for solving the problem. A list of \lqq invented" formulas for the integral is as follows: 

$\int_0^1 (x^2-y, y^2+x) ds$  $= (\int_0^1 x^2-y ds, \int_0^1 y^2+x ds)$;  

$\int_0^1 \vec F ds+\int_1^3 \vec F ds$; 

$\int_0^2\int_1^2 (1+x^2-y)dydx+\int_2^3\int_0^2 (-x+3)dydx=$\dots $=7+17=24$;

$\int^0_{1+x^2} (x^2-y, y^2+x)$,  $\vec f(x,y)=(x^3/3-y^2/2, y^3/3+x^2/2)$ from $C$ to $B$ to $A$. 

$\int_{-x+3}^{1+x^2}(x^2-y,y^2+x)$ and stops here.

$\int_0^3 (t^2-(1+t^2), (1+t^2)^3+t)\cdot\vec\nabla F(x,y)=$ (this comes after a correct parametrization of the parabola path, but no parametrization of the line path, and correct determination of $\nabla F(x,y)$). The incorrect answer obtained was 71187/70

$\int_c \vec F ds=\int_a^b (c(t)^2-c(t), c(t)^2+c(t))\cdot (c\pr (t)) dt$
These answers are a good indication that more work needs to be done in the area of vector field integral, certainly more examples to illustrate differences in integration of a vector field vs. scalar function. 

One student has the following incorrect solution: 
He draws the graph fairly correctly, and seems to parametrize the portions accordingly, but instead of evaluating the integral in the two portions, he  \lqq adds" the portions of the path as follows: $C_1=y=1+t^2+C_2=y=-(t-3)$, $c(t)=C_1+C_2=2y=t^2-t+4$, writes $c(t)=(t^2-t+4)/2$ and $c\pr(t)=(2t-1)/2$ and then $F(c)=(t^2-(t^2-t+4)/2, (-(t^2-t+4)/2)^2+t) $ and $F(c)\cdot c\pr=
(t^2-(t^2-t+4)/2, (-(t^2-t+4)/2)^2+t) \cdot (2t-1)/2$ and then $(2t-1)/2\cdot \int_{-1}^1 (t^2-(t^2-t+4)/2, ((t^2-t+-4)/2)^2+t)$.

 Most students answered correctly the orientation question, but there were a couple of incorrect answers as follows:

 \lqq The value wouldn't change going the other direction since the line didn't change shape or orientation."

"The value of the integral of the orientation of curve $\gamma$ went from C to B to A does not change as if it went from A to B to C." 

Here it would be worthwhile to spend extra time time with the one-dimensional case and why the sign of the integral changes when the orientation of the interval of integration is reversed. 

There are also some semi-correct answers, as follows: 

\lqq If the orientation was from C to B to A the value of the integration would be negative" (I corrected this to \lqq opposite").

This comes from the universal error (colloquialism) that students make, namely identification of \lqq opposite" with \lqq negative" which is not an innocent error, but is rather a linguistic sloppiness impeding students' understanding of concepts,  and an effort should be made to correct it.\footnote   {I have noted elsewhere that this error prevents students from adopting the definition of $|x|=$ $x$, if $x$ is positive and $-x$, if $x$ is negative, for $|x|$ is positive and $-x$ (\lqq negative x")\lqq is negative." }
\ju
\ju


{\bf Problem 5} The function $f(x,y)$ is defined as follows on the domain $[0,1]\times[0,1]$:

$$
f(x,y)=
\cases 
\frac{x}{\arctan (1/y^2)},  &\text{if}\quad  y\neq 0 \\
2x\pi, &\text{if}\quad y=0 \cr
 \endcases$$

Find all the points where the function is continuous (or, complementarily, where it is not continuous) in its domain. 

\sju

Here students have to display some sophistication, involving the definition of a continuous function as well as the knowledge of some trigonometry/algebra, such as finding when $\arctan$  is 0, knowledge of how to do some limits, etc. I have spent some classroom time explaining  that what they mostly (and erroneously) consider to be continuity (namely that the graph can be drawn without lifting the pencil off the paper), is not very useful in functions of several variables. I have emphasized on numerous occasions that the definition of continuity of $f$ at a point $P$ is that $\lim_{x\to P} f(x)=f(P)$. Although this is not such an overly complicated definition students do not readily accept it, for the simple reason that it is not the instant product they like -- finding limits may be an involved process and may look as a roundabout way to establishing what they like to think of as the intuitive concept of continuity. 

Certainly it would be nice to check that the function is defined in the proclaimed domain, by noting that $\arctan(1/y^2)$ never turns into $0$. The first stage of the proof would be to note that for points that have non-zero $y$-coordinate, the first part of the definition of $f$ is the quotient of two continuous functions. Students are expected to say that $x$ is continuous as a polynomial and that a rational  function with non-zero denominator is continuous as a composition of two continuous functions $\arctan$ and the rational function $1/y^2$. This argument so far would have earned them a good partial credit. Then they need to check continuity at points of the form $(a,0)$ i.e. they need to check that $\lim_{(x,y)\to(a,0)}f(x,y)=f(a,0)=2a\pi$. Students would have to find the  $\lim_{(x,y)\to(a,0)} x/\arctan(1/y^2)=a/(\pi/2)=2a/\pi$. This means, that in order to have continuity at such points, the equality $2a\pi=2a/\pi$ must be satisfied and that happens only for $a=0$. Hence the conclusion is that $f$ is continuous at all points $(x,y)$ for which $y\neq 0$, also at $(0,0)$ and discontinuous at all other points, namely all $(a,0)$, where $a\neq 0$. 

Here is how one student (the weaker of the two high school seniors who took the class) attempts to solve the problem:   She first rewrites the function in \lqq her own" way: $\tan^{-1}(1/y^2)$ and then further \lqq simplifies" the function into $x=\tan(1/y^2)\cdot y$, where the last $y$ stands for the function, namely $f(x,y)$.  My experience helps me to foresee this kind of error, so I  warn the students in classes that it is safer to use $\arctan$, rather than $\tan^{-1}$ lest they mistake it for the reciprocal, rather than the inverse function. However one has to spend additional time in clearing up this notational confusion.
Almost all widely used calculus textbooks  (see e.g. Rogawski (2008) or Stewart (2012) unwittingly or purposefully use the inverse function notation, lest they appear not to be sufficiently \lqq modern." It is as if we should be writing $square^{-1}$ to denote the square root function, just so we can show it is inverse of the basic quadratic function. The same bad habit is seen with calculator keys that almost exclusively use the inverse function notation. This notation is certainly useful in general discussions about inverse functions, but is a bad notation when names for inverse functions exist (for instance the $arc$ functions such as $\arctan$), which in their names have interpretation, if not the meaning of a particular inverse function. It has been noted many times before that good notation is golden and in this case it would result at least in students not confusing the inverse function with the reciprocal of the function.\footnote   {Most students are genuinely and thoroughly surprised to learn that $\tan^{-1}$ does not denote the reciprocal of $\tan$. It is an interesting functional equation to find all invertible functions $f$ whose inverse equals its reciprocal.  } 
To follow, the same student further draws that part of the square root function graph (in the 2D Cartesian coordinate system) that is in the unit square and claims it is the graph of $x/\arctan(1/y^2)$ -- all graphs for her are in two dimensions. 
Then continues on to say \lqq derivative of " $f\pr(x,y) =$ substitution $u=\tan^{-1}(1/y^2)$... $d 2x\pi/dx =2\pi$!! [sic]. Then claims that  the function is not continuous for $x=1$ and concludes: \lqq The graph $x/\arctan(1/y^2)$ lies in the positive $(x,y)$ plane and the function is discontinuous at $x=1$."
Yet again, more work must be done  simply on functions of more than one variable, in what their domains and codomains are, as well as their graphs, for, even if students have better understanding of functions than the action conception, without firm grasp of domains and codomains or ranges of participating functions would not lead them to correct solution.  

Another student evaluates the function at the four corners of the domain and concludes that \lqq the function $F(x,y)$ is continuous throughout its domain." This comes from students' confusion with procedures that were used in class when finding absolute extrema where points on the boundary of the region (or  at the corners of polygonal domains) were crucial. This is confirmed by another student finding the partials" $f_x, f_y$ (that seem to be incorrectly done  \lqq partial" anti-derivatives), which he then sets equal to zero and solves the resulting equations.  Unquestionably,  these confusions point to insufficient time that students had in digesting the material, be it objectively (one semester for all the advanced material) and/or subjectively -- insufficient work on the students' part).

One student correctly computes  the function values at the corners, 
then writes 
$$
f_x=
\cases 
 \frac{- x^2}{2\tan^{-1} (y^2)-\pi},  & \\
x^2\pi & \cr
\endcases$$

$$
f_y=
\cases
 -2  \frac{1}{2\tan^{-1} (y^2)-\pi},  & \\
2x\pi y & \cr
\endcases$$

He then solves when the parts of these \lqq partials" are equal to 0, conveniently getting $x=0$ as a solution to the first system and not getting anything (\lqq false") for the solution of the first part of the second expression; however, from the second part of the second system he gets $x\pi y=x,y=0$ and says \lqq these points are on the boundaries." This is a confusion with the problem of finding the extrema of functions;  the conclusion is that \lqq the function is continuous on the points (1,1) and (1,0) and not continuous at (0,0) and (0,1)."

One student attempts to find the double integral, over the square, of the first defining part of the function and, after some incorrect substitutions, etc., arrives at an incorrect value of the integral. Then he states, as a matter of fact that:   \lqq It is continuous because it is defined everywhere and limits exist."  Double integrals are frequent among the incorrect solutions here. 

Some students know the correct definition of continuity, but cannot find the limit of the function at a general point, rather they try to find it in a corner point or two, or invent their own limit rules. Thus one states:   \lqq$2x\pi$ is not dependent on $y$ and is a continuous function that will be well-defined within the domain" and then 
$$\lqq \lim_{y\to 0} \frac{x}{\arctan(1/y^2)}=\lim_{y\to 0} f(x,y)\pr=\frac{1}{\lim_{y\to 0}\text{arcsec}(1/y^2)}."$$

One student discusses continuity of the parts separately: \lqq For $f(x,y)=2x\pi$ this is continuous for all points work from the domain $x$= all real numbers. 
For $x/\arctan(1/y^2)$  all points work except for 0 from the domain. Because only way $x/\arctan(1/y^2)$ is not continuous is if the denominator is equal to zero. Only value it will equal zero is if $y=0$. So all points on the function $f(x,y)$ is continuous except when $y=0$."
This student did receive partial credit, but here comes again the problem of functions, in that these piecewise definitions of functions are often seen as two separate and independent functions, partly because it is more convenient to do so. Breidenbach, Dubinsky, Havks, Nichols (1992) use \lqq split" functions to test whether students' understanding of functions is the process one. This is something that needs to be addressed in pre-calculus or the beginning algebra courses. 

One student states: \lqq ...  a function is not continuous where it is undefined," then goes on to find correctly that the denominator in the fraction part is never $0$ and concludes that there are no points of discontinuity. 
This statement is an innocent regurgitation of classroom and textbook folklore where continuity of a function at a point is still considered even though there is no function (i.e. the function  is not defined there) -- the point addressed in  Dimitric (2004); see furthermore, Rogawski (2008) and Stewart (2012) where the instructing examples and the exercises on \lqq continuity"  are almost exclusively such that a student finds zeros of the denominators of rational functions and pronounces that the \lqq functions" are not continuous (since they are not defined) at those points. This is propagated as well into online homework systems such as Webassign, etc. 
\ju
\ju


{\bf Problem 6}  Find the volume of the region bounded by the hyperbolic cylinders $xy=1, xy=9, xz=4, xz=36, yz=25, yz=49$.
Do only the part in the first octant ($x,y,z\geq 0$). [Hint: Make an appropriate change of variables and use the change of variables theorem to compute it. ]

\sju

Setting up the problem and the change of variable theorem details are important and then the integral computation (with the Jacobian)... I have directed students to use change of variables to check whether they have learned this very process.  Equal points were allotted for substitution, for theorem conditions and for computation of the integral/Jacobian. 

Since each of the $x, y, z$ combinations repeat twice in a circular fashion, the substitutions seem to be imposing themselves: $u=xy, v=xz, w=yz$ with boundaries for $u|_1^9$, $v|_4^{36}$, $w|_{25}^{49}$ (*). Multiplying all the substitution equalities leads to a  helpful relationship: $xyz=\sqrt{uvw}$ (**). The map $\psi(x,y,z)=(xy,xz,yz)=(u(x,y,z), v(x,y,z), w(x,y,z))$ transforms a \lqq curvy, hyperbolic" $x,y,z$ brick $\Cal V$ into a perpendicular $u,v,w$ brick $\Cal V_0\,\,$ $[1,9]\times[4,36]\times[25,49]$. For $\psi$ to be a good substitution map, we need to check that it is one-to-one. That follows from the fact that equalities $x_1y_1=x_2y_2, x_1z_1=x_2z_2, y_1z_1=y_2z_2$ imply equalities $x_1=x_2, y_1=y_2, z_1=z_2$. Thus existence of the inverse function $\psi^{-1}=\phi(u,v,w)=$ 
$= (x(u,v,w), y(u,v,w), z(u,v, w))$ is guaranteed. In fact, solving the substitution equalities (*) gives $\phi(u,v,w)=(\sqrt{uv/w}, \sqrt{uw/v}, \sqrt{vw/u})$, which is also one-to-one. The Jacobian is then computed to be $\partial(u,v,w)/\partial(x,y,z)=-2xyz=-2\sqrt{uvw}$ ((**) was used). Both $\psi, \phi$ are $C^1$ maps, since all the partials exist and they are continuous. Thus conditions for the change of variable formula are satisfied and we find the volume as follows:
$$\iiint_\Cal V dx dy dz=\iiint_{\Cal V_0} |-2\sqrt{uvw}|^{-1} du dv dw$$
$$=\frac{1}{2}\int_1^9u^{-1/2}du\int_4^{36} v^{-1/2}dv\int_{25}^{49} w^{-1/2}dw=64.$$

There are a number of tripping points here, from routine to more subtle ones. Not knowing/remembering  how to compute the Jacobian (a fairly intuitive formula) was a rare miss among the routine ones. The change of variables formula (a substitution formula) is computed by 
$$ \iiint_\Cal V f(x,y,z) dx dy dz=$$
$$\iint_{\Cal V_0} f(x(u,v,w),y(u,v,w),z(u,v,w))|Jac(\phi)|du dv dw
$$
where $\phi:\Cal V_0\lra \Cal V$ is a $C^1$ and one-to-one map on the interior of $\Cal V_0$ and $f$ is a continuous integrand; here the Jacobian $Jac(\phi)=\prl(x,y,z)/\prl(u,v,w)$ (see e.g.[ Rogawski, 2008]). 
This is somewhat a convoluted formula, for it seems to work \lqq backwards" in that we do not know $\phi$, rather we know an easier $\phi^{-1}=\psi$ that is usually suggested by the given shapes of regions of integration. 
The instructor should spend time to justify the conditions and show that the formula does not work if conditions are not satisfied. He must further spend time in explaining the built in deformation of volume and thus justify the use of the Jacobian multiplier. The shortcuts students make are often utilitarian in nature and one of them is in not bothering with the Jacobian. Thus one student says: \lqq  If Jacobian was not needed then $\iiint dV=\int_1^9\int_4^{36}\int_{25}^{49} dV= \int_1^9\int_4^{36}\int_{25}^{49} dz dy dx$,  and so on and gets the volume 6144." Another student proceeds to do the same \lqq formula" without even mentioning the Jacobian. 

A technicality should also be explained, namely that the Jacobians of inverse functions are reciprocal to each other; some students had missed this fact. Emphasis on this however should be made already in the first calculus courses where derivatives of inverse functions are discussed.  

One of the students who missed including the Jacobian in the formula (but had the correct substitution) went further in simply calculating the volume of the straight box $\Cal V_0$ by multiplying its dimensions (which she found correctly from the substitution). 

Another student seems to know the general procedure and  he correctly chooses the substitution $u=xy, v=xz, w=yz$, but when he tries to solve for $x,y,z$ he is not successful. He puts a big \lqq ?" and writes \lqq no idea. If I knew the variables, I can figure it out, but I can't come up with the variables." Thus solving these algebraic systems is important and it is almost never done in lower division algebra classes, or is mentioned only as an afterthought. 

One student starts with $\int_{25}^{49}\int_4^{36}\int_1^9 r dr d\theta dz$ and proceeds with $r^2|_1^9$ etc to get something like 30720 at the end. Here, there is confounding with specialized procedures and a number of exercises done on the topic of special substitutions, such as cylindrical and spherical coordinates (with their own Jacobians); subtleties are missing though. This again means more exercise time should be devoted both in class and by the students so as to definitely discern the confusing notions and make them clearer in students' minds. 

One problem (not widespread) was not quite knowing the integral formula for the volume and it is tempting for some students to put $xyz$ or $uvw$ as the integrand in computing the volume (\lqq length times width times height"). 

There was one \lqq silent solution" (namely, no explanation, just a formula or two)  by a student who writes: 
$$ \int_1^3\int_5^7\int_2^6 z dz dy dx 
$$
and goes on to integrate this convenient integral with separated variables to get 64 as the solution. One can attribute this solution to a sharp talent to see through promptly, but, unfortunately this was not the case here. It is not an infrequent luck at numerology that almost with no exception accompanies problems that have suspiciously nice, round numbers (that, in this case also seem to be squares too). One can try to justify this solution by a substitution $u^2=xy, w^2=xz, v^2=yz$;  however since the one-to-one requirement is needed, the substitution is $u=\sqrt{xy}, w=\sqrt{xz}, v=\sqrt{yz}$, to get the correct limits for $u,v,w$.  When the Jacobian is found: $\prl(u,v,w)/\prl(x,y,z)=1/4$ and its reciprocal is 4, which will then give the correct result after evaluating
$$ \int_1^3\int_5^7\int_2^6 Jac\,\, dw dv du 
$$
None of this is to be found in the student's solution who had a lucky confounding event with a formula for finding volume under a surface $z$ (which would be a double integral). Two more students start taking square roots of the given numbers, for the temptation is too great not to do that when perfect squares of integers are given. 

It would help greatly teaching (and learning)the  integration method of substitution in several variables, if this method were thoroughly taught in the case of one variable (usually in calculus I). The practice of applying the method of substitution in one variable almost exclusively omits checking that the substitutions are of the right kind (namely that they are $C^1$ and one-to-one, etc.). Also, when doing substitutions, students retain the old limits since they want to go back to the old variable, once they eventually solve the indefinite integral via the substituting variable. Thus an opportunity is missed to introduce (implicitly) the \lqq Jacobian" and all the other features that would become rather useful when the case of several variables is taken up. I cannot overemphasize this point. 
\ju
\ju

%

{\bf Problem 7.} Find the volume of the smaller region bounded by the sphere $x^2+y^2+z^2=a^2$ ($a$ is a given constant) and the plane $z=b$, where constant $b$ is such that $a>b>0$, using either cylindrical or spherical coordinates.

\sju

This exercise could be worked out using elementary calculations, as long as the student knows formulas for volumes of balls, cones, etc. The essence of the computation is the method of cross sections, which is often taught in Calculus I courses. In the context of Calculus III, the problem is meant to check whether students have the skills to work out the volume via triple integrals, using specific substitutions, namely cylindrical or spherical coordinates. I specified these two methods perhaps to aid students in deciding how to start to solve, in view of the fact that either coordinate substitution would require about the same amount of work to find the volume. One student thought that one method was more advantageous than the other, so he was pondering which one to choose before he proceeded. Thus he thinks aloud: \lqq Region is bounded by a sphere so we will use spherical coordinates. The Jacobian or multiplier in spherical coordinates is $\rho^2\sin\rho$ [sic]" and then he lists spherical coordinates. But then crosses it all out and continues:  \lqq Actually since we are bounded by $z=b$, we will use cylindrical coordinates"... and lists polar coordinates for $x, y$, but not $z$ and writes that Jacobian $= r$. Unfortunately he does not find the limits of integration correctly. 

Setting up the volume integral earned one third of the assigned points: For cylindrical coordinates we have $x=r\cos\theta, y=r\sin\theta, z=z$ and then $V=\iiint_\Cal V dx dy dz=\iiint_{\Cal V_0}r dr d\theta dz$ and for the spherical coordinates $x=\rho\cos\theta\sin\phi, y=\rho\sin\theta\sin\phi, z=\rho\cos\phi$ one has $V=\iiint_{\Cal V_0} \rho^2\sin\phi d\rho d\theta d\phi$. These are substitution formulas with Jacobians already incorporated for these specific substitutions. Textbooks usually discuss them separately and my lecture time is also partly devoted to the same. If students cannot remember these formulas they have to use the general substitution procedures and compute the Jacobians, but they certainly have to know the spherical and cylindrical substitutions. 

The next stage is to determine $\Cal V_0$, namely to determine integration limits for the new variables. For this,  well-drawn geometric diagrams are useful, and simple appeal to the Pythagoras' theorem will give limits for cylindrical coordinates: $0\leq\theta\leq 2\pi$, $0\leq r\leq \sqrt{a^2-b^2}$ and $b\leq z\leq \sqrt{a^2-r^2}$. For spherical coordinates, the limits are as follows: $0\leq\theta\leq 2\pi$, $0\leq\phi\leq\arccos(b/a)$, and $a\leq \rho\leq b/\cos\phi$. Determining these limits was worth 1/3 of total points and so was the final part of actually computing the integrals (and thus the volume):
In cylindrical form
$$V=\int_0^{2\pi} d\theta\int_0^{\sqrt{a^2-b^2}} rdr\int_b^{\sqrt{a^2-r^2}}dz=$$
$$=2\pi\int_0^{\sqrt{a^2-b^2}} r dr (\sqrt{a^2-r^2}-b)=$$
$$=2\pi\left(\int_0^{\sqrt{a^2-b^2}} r\sqrt{a^2-r^2} dr - b\int_0^{\sqrt{a^2-b^2}}r dr\right).
$$
The first integral in parentheses is done by substitution $t=a^2-r^2$ to get the volume $V=\pi(b^3+2a^3-3ba^2)/3$.
Similarly, in spherical coordinates one has:
$$ \int_0^{2\pi} d\theta \int_0^{\arccos(b/a)} d\phi\int_{b/\cos\phi}^a \rho^2\sin\phi\, d\rho=$$
$$=2\pi\int_0^{\arccos(b/a)}\left( (a^3/3)\sin\phi-(b^3\sin\phi)/(3\cos^3\phi)\right) d\phi
$$
The second part of the last integral is done by substitution $t=\cos\phi$ and the volume again is $V=\pi(a-b)^2(2a+b)/3$. 

Here the variety of errors occur at every potential trouble spot. The instructor should spend more time in elaborating on the simple volume formula $V=\iiint_\Cal V dx dy dz$ and point out that the integrand is $1$, because the volume is already encoded in the region of integration $\Cal V$ whose volume is being calculated. Here students use all sorts of expressions for the integrand in this volume formula, most frequently  $x^2+y^2+z^2$ or $\sqrt{a^2-x^2-y^2}$. Other \lqq integrands" are invented by substituting either cylindrical or spherical coordinates into the sphere equation so that the integrands turn out to be $r^2+z^2$, $\sqrt{a^2-r^2}$, $\rho^2$ or  $\rho^2-a^2$. These errors come partly from students confounding procedures they used to find areas under curves (integrand is the function formula), or volumes under surfaces (integrand is the surface formula with the double integral as a device). Thus these should be contrasted with the volume formula when starting work with triple integrals.

With a number of errors along the way, students struggle through, coming up with some answers that are more often than not polynomials in $a, b$ (including numerical constants). However,  these polynomials range from being non-homogeneous, to being homogeneous of degrees 2, 4, 5 and so on. Teaching students basics of dimensional analysis should be done early on in (pre)algebra and introductory science courses, so that students are comfortable with various dimensions or units for that matter. It would then go a long way in recognizing problems with their \lqq solutions" simply based on the incorrect dimension they come up with for, in this case the volume. 

Determining limits of integration is definitely a hurdle and while no problem is found (with a couple of exceptions) in finding limits $0\leq\theta\leq 2\pi$, since most examples done in the textbook and classroom have those exact limits, determining limits for other variables is hard for many, not the least because students fail to make a good geometric drawing, but also because they are not yet good at understanding that limits of one variable need not be constant, i.e. they may depend on (an)other variable(s). Thus, in students' papers, $r$ is found to have limits from $0$ to $a$, or $b$ to $a-r$, or $0$ to $a-b$ and $z$ has limits  $0$ to $a$, or $0$ to $b$, or $b$ to $a-r$ or $0$ to $\sqrt{a^2-r}$, $b$ to $a^2+b^2$, or $a-b$ to $a$
$\phi$ varies from $a$ to $b$, $0$ to $\pi$, $0$ to $2\pi$.

The cursory role of 3D geometry in US school system is well-known. Even the best  of our students competing at math Olympiads need to be given extra preparation on this aspect of geometry. With that in mind, an instructor has to spend some time in three dimensional calculus drawing 3D pictures and manipulating them, through finding relationships, using similarities of triangles, and all. Computer aid in 3D imaging is pretty, but perhaps not very helpful, for students still would not be able to draw figures and subfigures, find the limits and relationships among various elements in their 3D figures, unless they made repeated attempts at drawing them by hand. With slight exaggeration, one should require simple art classes as a requirement for taking a tad more advanced math class. 

One student has a definite idea to solve the problem as a double integral, namely to play with volumes under the spherical cap and the flat cap $z=b$, all using cylindrical coordinates.  But she  says \lqq $x^2+y^2+z^2=a^2$ projection onto {\it x-y} plane $=x^2+y^2=a^2$," rather than projecting the circle at level $z=b$. 
\ju

\centerline{\bf A brief lesson plan}
\ju

One of the many uses of our analysis is in making lesson plans. 
Let us assume that the lesson's topic is continuity, one of the most important and most neglected topics in calculus. We have feedback here attached to Problem 5 and we can incorporate it into the lesson plan as follows:

Make sure that, by this time,  students are well-familiar with the exact notion of function, as an ordered triple -- domain, codomain and the assignment rule(s). This should be done at early pre-calculus,  but in any case it is a must to review/teach at the beginning of the course. Give students some \lqq trick questions" such as ask them to decide which of the following are functions and which are not: 
$f:A\lra B$,   a) $f(x)=\sqrt x$, $A=\RR$, $B=\RR$, b) $f(x)=\sqrt x$, $A=\RR^+$, $B=\RR$, c) $f(x)=-\sqrt x$, $A=\RR^+$, $B=\RR^+$, d) $f(x)=-\sqrt x$, $A=\RR^+$, $B=\RR^-$, e) $f(x)=\sqrt x$, $A=\RR^+$, $B=\RR^-$. If complex numbers are part of the course, tailor this accordingly.  

By this time enough exercises must have been done for students to feel comfortable with the notion of a limit and the routine of finding various limits, as well as basic operations with limits. This should already be done in the first calculus course. Test them on this with the following limits:  
$$ \lim_{x\to 0} \frac{\sin x}{x}, \lim_{x\to 0} \frac{\cos x}{x}, \lim_{x\to\infty} \frac{\sin x}{x}.
$$
It is likely that the students will indiscriminately invoke \lqq L'Hospital's rule;" remind them what conditions must be satisfied with this rule. 

Now, use some limits, in higher dimensions, preferably those you will use to check whether functions are continuous or not. 
First some that separate variables: 
$$ \lim_{(x,y)\to (0,\pi/2)}\frac{\sin x\cos y}{xy}, 
$$
And some that are not obviously variable separating $\lim_{(x,y)\to (a,b)}[xy]$ (the integer part function). Discussion of cases when $ab$ is and is not an integer may be instructive. 

Do also some piecewise defined function examples, for instance the example used in Problem no. 5. Or the following one: 

$$
f(x,y)=
\cases 
 \frac{e^x-1}{x}\frac{\sin y}{y},  &  \text{ if } x\neq 0 \text{ and } y\neq 0, \\
 \frac{\sin y}{y},  &  \text{ if } x=0 \text{ and } y\neq 0, \\
 \frac{e^x-1}{x},  &  \text{ if } x\neq 0 \text{ and } y=0, \\
 1,  &  \text{ if } x=0 \text{ and } y=0. \\
\endcases$$

Define continuity of function $f$  at a point $P\in A$ of the domain of the function $f:A\lra B$ if $\lim_{x\to P}f(x)=f(P)$. Emphasize that a question whether a \lqq function" is continuous at a point outside of its domain is meaningless, since there is no function at such a point. Define what is meant to say that $f$ is discontinuous at a point.  Define (global) continuity of $f$ on all of its domain. 

Now do first some simple examples of one variable, then pick up the aforementioned examples from limit exercises and add new ones. Make sure the examples are graduated, from simple to more complicated. 

If possible, visualize the functions and their limits, either by hand, or using graphing software. 

\vfill
\eject

\centerline{\bf Conclusions}
\sju

The problems with understanding functions exist from early mathematics education and are dragged along throughout, including higher level or even graduate level courses. I believe that the foundation for these problems begin before functions are introduced, at the level of basic notions such as sets and ordered pairs (or triples, etc). Graphing of functions is introduced simultaneously with the Cartesian plane, without much ado about abstract formation of the Cartesian product of (two) sets. Yet, it took many generations in human development before the idea of Cartesian product matured. Not only should considerable time be devoted to \lqq playing" with coordinates, but a physical drawings\footnote   {I like slightly exaggerating and saying that \lqq what does not go through your hand does not get to your head..."} of coordinates (and at least simple objects) in two and three dimensional coordinate systems should be practiced thoroughly until the students feel this as a \lqq second nature." Once this ground work is made for the dimensions higher than one, one can introduce a name for a subset of a Cartesian product $A\times B$ -- call a subset $\rho\se A\times B$ a {\it relation from $A$ to $B$}, that can be denoted by $\rho:A\lra B$.\footnote   {My own experience and feedback from several educational systems show that even very young students are not intimidated by this abstraction level, which may run contrary to prejudices in existence regarding mathematics education.} Some simple operations and drawings with relations should be introduced as a practice. 

Then, time comes for a definition of a function $f:A\lra B$ as a relation with properties that, for every $a\in A$ there is only one $b\in B$ with $(a,b)\in f$ or $b=f(a)$. $A$ should be called the {\it domain} of $f$ and $B$ the codomain and the range would be introduced as the set of all $b\in B$ such that $(a,b)\in f$, for all $a\in A$. An important point should be made, namely that the function consists of these three parts $f:A\lra B$ and that each of them as important as any other in their own right and that two functions $f:A\lra B$ and $g:C\lra D$ are equal if and only if all three components are equal: $A=C, B=D, f=g$. A good number of exercises that vary some or all parts of a function and discussion of resulting differences should be made until students adopt this formalism to a high degree. 

Then the stage is set for metaphors.  A dynamic metaphor that I use is that of a processor (a \lqq grinding machine") as sketched in the following diagram: 

\sju
[Here a diagram of the processor needs to be inserted...]

One can represent this as making of bread sticks (placed in codomain $C$) after dough in domain $A$ is processed by pasta machine $f$.\footnote   {One can vary the description dependent on the fashionable language of the day and place -- veggie burgers in California or a sausage machine in Texas...} Here, it is very important to impart to students that a function $f:A\lra B$ is the whole assembly line consisting of $A$, $f$ and $B$, rather than just rely on a vague hope that the students would see functions either as action or process conception. Importance of domain may be imparted for instance by replacing domain with dough with domain with iron bars that cannot be processed by the processor, or codomain that is a capillary dish that cannot accommodate the processed dough... Varying all the parts of the function would show a student that each part is in fact equally important as a constituent of what is meant by a function. Composition may be explained in a similar manner, for a packing machine $g:C\lra D$ may be the follow up \lqq factory" that will work, only if the final product of $f$ is in the domain of $g$; thus $g$ may be able to process sticks only of certain length range, not of other dimensions, yet again underlining importance of domain and codomain (or range) -- the diagram has some sticks separated from the other in $B$, those in domain $C$ of $g$ and those that are not in $C$. 

There are other metaphors one can use, but a very important point to remember is that we cannot have correct metaphors if we do not give students correct definitions of notions we teach them and insist that they adopt these definitions, at least at the memorization level (first step in Bloom's taxonomy). Otherwise we would end up in finding metaphors, of metaphors, thus giving students an impression that the notions are negotiable moving targets. While mathematical formalism is not equivalent to what mathematics is, it is nonetheless an essential constituent of mathematics. It is most crucially visible in at least one practical development, namely computerization of human knowledge and experience which could not take place without strict mathematical formalism. Development of ideas is mostly a meandering process, not infrequently through bouncing off of the extreme views or practices. Thus an increased formalization of mathematics in the 20th century led to sometimes indiscriminate transplantation of formalism into teaching practice, which in turn led to overreaction in other direction, in denial of all the formalism where the notions of function and its properties was degraded to pre-Eulerian levels of understanding of what a function is. In a pivotal exposition of the pre-1930's history of the development of function notion, Luzin (1935; see also Yushkevich (1966, 1977)) highlights subtleties of the function concept inclusive with the then developing views on classification of functions, intuitionistic developments, etc. Bourbaki formalized the function concept in 1939 and then category theory brought about a great generality and the accompanying simplification of the concept. Adoption of these developments in school curriculum has been checkered at best and it would be an interesting topic of research to look into the reasons why function concept is not adopted and properly taught in schools (and colleges).  

If this huge leap of students remembering correct definition of function gets achieved, other constructs with functions would become much easier, one of them being definition of continuity. 

There are some drawbacks in analyzing or publishing an analysis of this kind. For instance, it is a fact that analyzing tests in detail is very time consuming, and  even tedious. In addition, there are instances where this powerful educational tool may be misused: Fearing negative feedback from students or administration alike, some instructors would attempt to avoid potential trouble spots (as found through analyses of this kind)  by considerably diluting the level of content of the courses both in the classroom as well as in levels of examinations. 

As with every powerful weapon, dangers of its misuse always exist. Still, the benefits here outweigh the dangers. Based on the work students have shown on tests one can extract a good amount of useful information that can be utilized as follows (beyond the obvious and customary assessment value of tests): 
\ju

\roster
\item "1." To find out about students' specific difficulties in learning specific topics. In this case we learn about problems students are likely to have with analytic geometry, linguistic and phonetic hurdles in adopting concepts, difficulties with \lqq higher" dimensions, differentiation and integration hurdles, (mis)understanding of the concept of function, and continuity, especially when mapping one higher dimension to another, in graphical visualization and (in)ability to use it intelligently. In (lack of) algebraic skills and ability to perform simple mental calculations. In problems in doing integration, drawing graphs in 2D, dimensional analysis, etc. 
\item "2." To learn about points of emphasis and amount of  repetition the instructor should devote to specific topics.
\item "3." To help instructors devise appropriate lesson plans, that would preempt usual difficulties, misconceptions and obstacles that exist before and during instruction on any specific topic. 
\item "4." To get a specific content for making diagnostic tests for any  particular class. 
\item "5." To learn about students' gaps in their knowledge in material they should have learned by the time they come to the specific subject (prerequisites). 
\item "6."  To get very concrete snapshots of progression of teaching and learning, in fairly real time (since midterm exams are part of the analysis).
\endroster
\sju
{\bf Topics for further exploration}
\sju
a) Compare performance of students in mathematics classes dependent on whether they had previously taken or not a drawing class, with perspective drawing in particular. A default expectation is that the students who took such a class would perform considerably better in multi-dimensional calculus or in other courses where more than one dimension is involved. A worthy cause would also be to quantify the difference in math abilities of the two groups. 

b) Compare performance of students in mathematics classes dependent on whether they had previously taken or not a mathematics culture class where considerable time is devoted to issues of coordinate systems (dimensions, ordered tuplets, relations and so on), variables, constants, sets and operations with them. My belief is that students who took such a class would be considerably ahead of those who did not take that kind of class.

c) Develop a calculus complex (consisting of three segments, say) in a way that earlier segments set up the stage for easier sailing into the subsequent segments, by way of generality of introducing concepts and by thoughtful process of choosing expandable metaphors. 
\ju
{\bf Acknowledgment:} The author would like to thank Professor Annie Selden for her useful feedback on the content and form of the paper, 
and Professor Ivko Dimitri\'c for their stylistic comments, as well as anonymous referees for their comments that led to a better presentation of the content of this paper. 
\ju


\centerline {\bf References}
\vskip 1truecm

Anderson, L.W. \& Krathwohl, D.R. \& Airasian, P.W. \& Cruikshank, K.A. \&
Mayer, R.E. \& Pintrich, P.R. \& Raths, J. \& Wittrock, M.C. (Eds.) (2001): A
taxonomy for learning, teaching, and assessing. A revision of Bloom's taxonomy
of educational objectives; Addison Wesley Longman, Inc.


Bloom, B.S. (Ed.) (1956). Taxonomy of Educational Objectives: The Classification
of Educational Goals, pp.201--207; Susan Fauer Company, Inc.

Breindenbach, Daniel \& Dubinsky, Ed \& Hawks, Julie \& Nichols, Devilyna (1992): Development of the process conception of function. {\it Educational studies in mathematics} {\bf 23}, 247--285.

Dimitric, R.M. (2001): Using less calcus in teaching calculus. A historical approach. {\it Mathematics Magazine} {\bf 74}, No.3 (June), 201-211.

 Dimitric, R.M. (2003): Components of Successful Education, {\it Tenth International Congress of Mathematics Education (ICME10)}, TSG 27, 2004; 
{\it The Teaching of Mathematics}, {\bf 6}(2003), No.2, 69--80 

Dimitric, R.M. (2004): Is $1/x$ continuous at 0?, arXiv: 1210.2939v1;mathGM] Preprint. 

Dimitric, R.M. (2009):  Student placement in calculus courses. A case study.
      {\it Eleventh International Congress of Mathematics Education (ICME11)}, TSG , 2009;  
      {\it The Teaching of Mathematics}, {\bf 12}(2009), No. 2,  83-102.

Dimitric, R.M. (2012a): How to make a course specific diagnostic test in the context of Bloom's taxonomy, In preparation.

Dimitric, R.M. (2012b): On complexity distance and its applications in education, In preparation.

Dubinski, E. \& Schoenfeld, A.H. \& Kaput, J.J. (Eds) (1994): {\it Research in Collegiate Mathematic Education, I}. Providence, RI: American Mathematical Society.

Dubinski, E. \& Harel, G. (1992): The nature of the process conception of function. In G. Harel \& E. Dubinsky (Eds.). The concept of function: Aspects of epistemology and pedagogy. {\it MAA Notes} {\bf 25}, 85--106. 

Horscovics, N. (1982): Problems related to the understanding of functions, in {\it Proceedings of the Workshop on Functions}, Organized by the Foundation for Curriculum Development, Enschede. 

Leinhardt, G. \& Zaslavsky, O. \& Stein, M. (1990): Functions, graphs and graphing: Tasks, learning and teaching, {\it Review of Educational Research} {\bf 60}(1), 1-64.

Luzin, N.N. (1935): Function, The Great Soviet Encyclopaedia {\bf 1959}, 314-334
Translated by Abe Shenitzer: Function: Part I {\it Amer. Math Monthly} {\bf 105}, No. 1 (Jan. 1998), 59--67 and No. 3 (Mar. 1998), 263--270). 

Oehrtman, Michael \& Carlson, Marylin \& Thomson, Patrick W. (2008): Foundational reasoning abilities that promote coherence in students' function understanding, In: Making the Connection: Research and Teaching in Undergraduate Mathematics, SIGMA and RUME, {\it MAA Notes}, Vol. 75, 27--41.

Piaget, J. \& Grize, J.-B. \& Szeminska, A. \& Bang, V. (1977): {\it Epistemology and Psychology of Functions} (J. castellanos and V. Anderson, trans.), Reidel, Dordrecht. 

Rogawski, J. (2008): {\it Calculus. Early Transcendentals}, 3rd printing. New York: W.H. Freeman and Company.

Selden, J. \& Mason, A. \& Selden, A. (1989): Can average calculus students solve nonroutine problems? {\it Journal of Mathematical Behavior} {\bf 8}, 45--50.

Selden, J. \& Selden, A. \& Mason, A. (1994): Even good calculus students can't solve nonroutine problems. IN: {\it Research issues in undergraduate mathematics learning} (Vol. 33, pp.19-26). Washington, D.C.: Mathematical Association of America.

Selden, Annie \& Selden, John \& Hauk, Shandy \& Mason, Alice(2010): 
Why Can't Calculus Students Access their Knowledge to Solve Nonroutine Problems? {\it CBMS Issues in Mathematics Education}, {\bf 8}, 128--153.

Stewart, J. (2012): {\it Single Variable Calculus. Early Transcendentals.} 7th ed. Australia: Brooks/Cole, Cengage Learning.

Sfard, A. (1987): Two conceptions of mathematical notions: Operational and structural, in J.C. Bergeron, N. Herscovics and C. Kieran (eds.), {\it Proceedings of the Eleventh International Conference for the Psychology of Mathematics Education, III}, Universit\'e,  Montreal, pp 162-169.

Tall, D. \& Vinner, S. (1981): Concept images and concept definitions in mathematics with particular reference to limits and continuity, {\it Educational Studies in Mathematics} {\bf 12}, 151--169. 

Thomson, Patrick W. (1994): Students, functions, and the undergradute curriculum. In Dubinsky et al (Issues in Mathematics Education Vol 4, pp 21-24).

Vinner, S. (1989): Images and definition for the concept of function,
{\it Journal for Research in Mathematics Education} {\bf 20}, 356--366.

Vinner, S. (1992): The function concept as a prototype for problems in mathematics learning, in Dubinsky \& Harel (1992), 195--214.

Vinner, S. \& Dreyfus, T. (1989): Images and definitions for the concept of function, {\it Journal of Research in Mathematics Education} {\bf 20}. 356--366. 

Yushkevich, A.P. (1966): The concept of function up to the middle of the nineteenth century,  {\it Istor.-Mat. Issled.} , {\bf 17},  123--150.

Yushkevich, A.P. (1977): The concept of function up to the middle of the nineteenth century,  {\it Arch. History of Exact Sci.} , {\bf 16} 37--85.

\vfill
\eject

\centerline{\bf APPENDIX: CALC III DIAGNOSTIC TEST}
\ju

What was the last math course you took, when and where?

Career goal (what is the reason you are taking this course):

Please do the following problems in one hour, without using a
calculator. These problems are a quick sample of calculus and algebra
you will need for Calculus III. I want to see how fresh you are with
the material and advise you accordingly.


\ju
1) Sketch the curves $y=-2x-1$ and $y=x^2-1$ and find
algebraically the points of intersection of the two lines.
 

2) Define when a function $f(x)$ is continuous at a given point $a$.

3) Solve the system of equations:
$$x+y=1,\quad x^3+y^3=2.\hskip 8truecm$$

4) Find the limit 
$$\lim_{x\to 0}\frac{\sin{3x}}{5x}=\hskip 8truecm$$

5) What is the set of all numbers $x$ such that $e^{-x}>3$?

6) Draw the graphs of the functions $e^{-x}$ and $ln(-x)$.

7) What are all the $x$ such that $\sin x>\cos x$?

8) What does the equation $(x+1)^2+(y-1)^2=1$ represent geometrically?

9) What is the volume of a right circular cone with radius of base $r$ and height $h$?

10) Find
$$\lim_{x\to\infty}\frac{2x^2+3x^7+1000}{2x^7+3x^2+1001}=\hskip 6truecm$$

11) Write the equation of the tangent line to the graph of the function
$$f(x)=\frac{x}{1+x^2},$$ passing through the point (1,1/2).

12) Find the first derivative of the function $f(x)=(-2x+1)^{100}$

13) Calculate $\sin^2 70^\circ+\cos^2 70^\circ=\hskip 8truecm$

14) True or false: $\sqrt{a^2+b^2}=a+b ?\hskip 8truecm$

15) True or false: 
$$\frac{a}{b+c}=\frac{a}{b}+\frac{a}{c} ?\hskip 8truecm$$

16) Calculate the integral: 
$$\int xe^xdx=\hskip 8truecm$$

17) If $f(x)=ae^{2x^2}$, what is $f(2a)$ ? ($a$ is a given constant).

18) 
Calculate the integral 
$$\int^1_0 \frac{x}{1+x^2}dx=\hskip 8truecm$$  

19) Calculate $\int^\infty_1 xe^{-x^2}dx=\hskip 8truecm$

20) Calculate 
$$\int_0^{\pi/4}\tan x\, dx=\hskip 8truecm$$

\enddocument